\documentclass[11pt]{amsart}
\usepackage{amsbsy,amsfonts}
\usepackage{latexsym,eucal,amsmath,amsthm}
\usepackage{amssymb}
\newcommand{\db}{\mathbb }

\newcommand{\dbR}{{\db R}}

\def\Im{\mbox{\rm Im}}

\def\X{{\bf X}}
\def\R{{\db R}}
\def\cal{\mathcal}

\theoremstyle{plain}
      \newtheorem{theorem}[subsection]{Theorem}
      
      \newtheorem{proposition}[subsection]{Proposition}
      \newtheorem{lemma}[subsection]{Lemma}

\theoremstyle{remark}
      \newtheorem{remark}[subsection]{Remark}

\theoremstyle{definition}
      \newtheorem{definition}[subsection]{Definition}

\begin{document}
\date{}
\title[Global well-posedness for Schr\"odinger with derivative]
{A refined global well-posedness result for
Schr\"odinger equations with derivative}

\author{J. Colliander}
\thanks{J.E.C. was supported in part by an N.S.F. Grant DMS 0100595.}
\address{\small University of Toronto}

\author{M. Keel}
\thanks{M.K. was supported in part by N.S.F. Grant DMS 9801558.}
\address{\small University of Minnesota, Minneapolis}

\author{G. Staffilani}
\thanks{G.S. was supported in part by N.S.F. Grant DMS  0100375 and
grants from Hewlett and Packard and the Sloan Foundation.}
\address{\small Brown University and Stanford University}

\author{H. Takaoka}
\address{\small Hokkaido University}
\thanks{H.T. was supported in part by J.S.P.S. Grant No. 13740087}

\author{T. Tao}
\thanks{T.T. is a Clay Prize Fellow and was supported in part by
a grant from the Packard Foundation.}
\address{\small University of California, Los Angeles}

\begin{abstract}
In this paper we prove that  the 1D Schr\"odinger equation with
derivative in the nonlinear term is  globally well-posed
in $H^{s}$, for $s>\frac12$ for data small in $L^{2}$.
To understand the
strength of this result one should recall that for $s<\frac12$
the Cauchy
problem is ill-posed, in the sense that uniform continuity with
respect to the initial data fails.
The result follows from the method of almost conserved energies,
an evolution  of the
``I-method'' used by the same authors
to obtain global well-posedness for $s>\frac23$.
The same argument can be used to prove that
any quintic nonlinear defocusing Schr\"odinger
equation on the line is globally well-posed for large data in
$H^{s}$, for $s>\frac12$.
\end{abstract}

\maketitle

\section{Introduction}

In this paper, using the method of almost conserved
energies,
we establish a sharp result  on global
well-posedness for
the derivative nonlinear Schr\"odinger IVP
\begin{equation}
\left\{ \begin{array}{l}
i\partial_tu+\partial_{x}^{2}u=i\lambda\partial_{x}(|u|^2u),\\
      u(x,0) = u_{0}(x),\hspace{1.5cm}x \in \dbR, \, t \in \dbR,
\end{array}\right.
\label{ivp}\end{equation}
where $\lambda \in \dbR$.

The first result of
this kind was obtained in the context of the KdV and the
modified KdV (mKdV) initial value problems (IVP) \cite{ckstt:2}, also using 
almost conserved energies.
Below we will discuss in more details the
``almost conservation method''
and its  relationship with the ``I-method'' which was  applied to \eqref{ivp}   
in \cite{ckstt:0} (see also \cite{keel:wavemap, keel:mkg, ckstt:1}).

  From the point of view of physics the equation in \eqref{ivp} is a model for
the propagation of circularly polarized Alfv\'en waves in magnetized
plasma with a constant magnetic field \cite{momt, m, ss}.

It is natural to impose the smallness condition
\begin{equation}\label{small-l2}
\| u_0 \|_{L^2} < \sqrt{\frac{2\pi}{|\lambda|}}
\end{equation}
on the initial data, as this will force the energy to be positive via
the sharp Gagliardo-Nirenberg inequality.  Note that the $L^2$ norm
is conserved by the evolution. In this paper, we prove the following
global well-posedness result:
\begin{theorem}\label{main}  The Cauchy problem \eqref{ivp} is
globally well-posed in $H^s$ for $s > \frac12$, assuming the smallness
condition \eqref{small-l2}.
\end{theorem}
We present here once again \cite{ckstt:0}
a summary of the well-posedness story for \eqref{ivp}.
Scattering and well-posedness for this  Cauchy problem has been
studied
by many authors \cite{hayashi, hayashi-ozawa, hayashi-ozawa2,
hayashi-ozawa4,
kaup-newell, ozawa,
ot, takaoka:dnls-local, tsutsumi-fukuda, tsutsumi-fukuda-2}. The best
local well-posedness result is due to
Takaoka~\cite{takaoka:dnls-local}, where a gauge transformation and
the Fourier restriction method are used to obtain local well-posedness
in  $H^{s}, \, \, s\geq \frac12$. 
In \cite{takaoka:dnls-global}, Takaoka
showed this result is sharp in the sense that the nonlinear evolution
$u(0) \mapsto u(t)$, thought of as a map from $H^s$ to $H^s$ for some fixed
$t$, fails to be $C^{3}$ or even uniformly $C^0$ in this topology, even when
$t$ is arbitrarily close to zero and the $H^s$ norm of the data is
small (see also Bourgain \cite{borg:measures} and Biagioni-Linares \cite{bl}).
Therefore,
we see that Theorem 1.1 is sharp, in the sense described above, 
except for the endpoint.

In \cite{ozawa}, global well-posedness is obtained for \eqref{ivp} in
$H^1$ assuming the smallness condition \eqref{small-l2}. The argument
there is based on two gauge
transformations performed in order to remove the derivative in the
nonlinear term and the conservation of the Hamiltonian.
This was improved by
Takaoka~\cite{takaoka:dnls-global}, who proved global well-posedness
in
$H^{s}$ for $s>\frac{32}{33}$ assuming \eqref{small-l2}. His method
of proof is based on the idea of Bourgain~\cite{borg:refinements,
borg:book} of estimating separately the evolution of low frequencies
and of high frequencies of the initial data. In \cite{ckstt:0}, we
used the ``I-method'' to push further the Sobolev
exponent for global well-posedness down to
$s>\frac23$. The main idea of the  ``I-method'' consists of
defining a modified $H^s$ norm permitting us to capture some
nonlinear cancellations in frequency space during the evolution
\eqref{ivp}. These cancellations allow us to prove that the modified
$H^s ({\mathbb{R}} )$ norm is nearly conserved in time, and an iteration
of the local result proves global well-posedness provided $s> \frac{2}{3}.$
In this paper, an algorithmic procedure, first developed in
the KdV context \cite{ckstt:2}, is applied to better capture the cancellations
in frequency space. Successive applications of the algorithm generate
higher-order-in-$u$ but lower-order-in-scaling corrections to the
modified $H^s$ norm. After one application of our algorithm, we
show that the modified $H^s$ norm {\it{with}} the generated
correction terms changes less in time than the modified $H^s$ norm
itself, so the first application of the algorithm
produces an {\it{almost conserved energy}}. The improvement obtained
allows us to iterate the local result and prove global well-posedness
in $H^s ( {\mathbb{R}} )$ provided $s> \frac{1}{2}$. In principle,
the algorithm may itself be iterated to generate a sequence of
almost conserved energies giving further insights into the dynamical
properties of \eqref{ivp}.
The end point
$s=\frac12$ is not obtained here. We speculate however that a
further  refinement of the `almost conservation method''
could be a possible way
to approach  this question.

We conclude this section with the following remark.

\begin{remark}
Consider the 1D  quintic nonlinear Schr\"odinger 
\begin{equation}\label{quintic}
    i\partial_{t}u=\partial_{xx}u +iau\bar u\partial_{x}u 
+ibu^{2}\partial_{x}\bar u+cu^{3}\bar u^{2},
\end{equation}
where $a, b$ and $c$ are fixed real numbers. If $(a+b)(3a-5b)/48+c/3<0$ the equation 
in \eqref{quintic} is defocussing and, as was remarked in \cite{ckstt:0}, 
the techniques used to prove Theorem \ref{main} apply here too and one 
can  prove global well-posed
for initial data in $H^{s}, s>\frac12$. Moreover, if $a=b=0$, 
we expect our method to give global well-posedness\footnote{Recall 
that in this case the initial value problem is locally well-posed in 
$H^{s}$, for $s\geq 0$, see \cite{cw}  and \cite{ts}.} even below $s=1/2$. 

We should point out that 
Clarkson and Cosgrove \cite{CC} (see also \cite{as}) proved that \eqref{quintic} fails  the Painlev\'e 
test for complete integrability when
$$c\ne \frac{1}{4}b(2b-a).$$
In particular this shows that our techniques, which do not depend on 
$a,b,c$,  do not rely on complete integrability.
\end{remark}

\section{Notation and Known Facts}

To prove Theorem \ref{main} we may assume
$\frac12< s \leq \frac23$, since for
$s > \frac23$ the result is contained in \cite{ozawa,
takaoka:dnls-global} and \cite{ckstt:0}.  Henceforth $\frac12< s
\leq \frac23$ shall be fixed.  Also, by rescaling $u$,
we may assume $\lambda = 1$.

We use $C$ to denote various constants depending on $s$; if $C$
depends on other quantities as well, this will be indicated by
explicit subscripting, e.g. $C_{\|u_0\|_2}$ will depend on both $s$
and $\|u_0\|_2$.  We use $A \lesssim B$ to denote an estimate of the
form $A \leq C B$,  and $A \sim B$ for $cB\leq A\leq C B$, where $c$ 
and $C$ are absolute constants. We also use $A\ll B$ if $A\leq 
\epsilon B$, where  $\epsilon$ is a very small absolute constant.  We use $a+$ and $a-$ to denote expressions of the
form $a+\varepsilon$ and $a-\varepsilon$, where $0 < \varepsilon \ll
1$ depends only on $s$.

We use $\|f\|_{p}$ to denote the $L^{p}(\dbR)$ norm, and $L^q_t
L^r_x$ to denote the mixed norm
$$ \| f\|_{L^q_t L^r_x} := \left(\int \|f(t)\|_r^q\ dt\right)^{1/q}$$
with the usual modifications when $q=\infty$.

We define the spatial Fourier transform of $f(x)$ by
$$ {\cal F}(f)(\xi):= \hat f(\xi)
:= \int_{\R} e^{-i x \xi} f(x)\ dx$$
and the spacetime Fourier transform $u(t,x)$ by
$$ \widetilde{{\cal F}}(u)(\tau, \xi) :=
\tilde u(\tau, \xi) := \int_{\R}\int_{\R} e^{-i (x \xi + t \tau)}
u(t,x)\ dt dx.$$
Note that the derivative $\partial_x$ is conjugated to multiplication
by $i\xi$ by the Fourier transform.

We shall also define $D_x$ to be the Fourier multiplier with symbol
$\langle \xi \rangle := 1 + |\xi|$.  We can then define the Sobolev
norms $H^s$ by
$$ \| f \|_{H^s} := \| D_x^s f\|_2 = \| \langle \xi \rangle^s \hat f
\|_{L^2_\xi}.$$
We also define the spaces $X^{s,b}(\R \times \R)$ (first introduced
in the contest of the Schr\"odinger equation in \cite{borg:xsb})
on $\R \times \R$ by
$$ \| u \|_{X^{s,b}(\R \times \R)} := \| \langle \xi \rangle^s
\langle \tau - |\xi|^2 \rangle^b \hat{u}(\xi,\tau)\|_{L^2_\tau L^2_\xi}.$$
We often abbreviate $\| u\|_{s,b}$ for $\| u \|_{X^{s,b}(\R \times
\R)}$.
For any time interval $I$, we define the restricted spaces $X^{s,b}(I
\times \R)$ by
$$ \| u \|_{X^{s,b}(I \times \R)} := \inf \{ \| U \|_{s,b}:
U|_{I \times \R} = u \}.$$
We shall take advantage of the Strichartz estimate,
(see e.g. \cite{borg:xsb})
\begin{equation}\label{strichartz-6}
\| u \|_{L^6_t L^6_x} \lesssim \| u \|_{0,\frac12+},
\end{equation}
which interpolates with the trivial estimate
\begin{equation}\label{strichartz-0}
\| u \|_{L^2_t L^2_x} \lesssim \| u \|_{0,0},
\end{equation}
to give
\begin{equation}\label{strichartz-p}
\| u \|_{L^p_t L^p_x} \lesssim \| u \|_{0,\alpha(p)},
\end{equation}
for any $p\in [2,6]$ and $\alpha(p)=\frac{(3+)(p-2)}{4p}$.
We also use
\begin{equation}\label{strichartz-2}
\| u \|_{L^\infty_t L^2_x} \lesssim \| u \|_{0,\frac12+},
\end{equation}
which together with  Sobolev embedding gives
\begin{equation}\label{strichartz-infty}
\| u \|_{L^\infty_t L^\infty_x} \lesssim \| u \|_{\frac12+,\frac12+}.
\end{equation}
The next lemma introduces two more estimates, that are probably less
known than the standard Strichartz estimates:
\begin{lemma}
For any $b>\frac12$ and any function $u$ for which the right hand
side is well defined, we have
\begin{equation}
     \|D_{x}^{\frac12}u\|_{L^{\infty}_{x}L^{2}_{t}}\lesssim
     \|u\|_{X^{0,b}}.
\label{se}\end{equation}

\noindent
(smoothing effect estimate).

For any $s>\frac12$ and $\rho\geq\frac14$ we  have
\begin{eqnarray}
  \label{mf1}\|u\|_{L^{2}_{x}L^{\infty}_{t}}
  &\lesssim& \|u\|_{X^{s,b}},\\
  \label{mf2}\|u\|_{L^{4}_{x}L^{\infty}_{t}}
  &\lesssim& \|u\|_{X^{\rho,b}},
  \end{eqnarray}

  \noindent
  (maximal function estimates).
\end{lemma}
\begin{proof}
The estimates \eqref{se}, \eqref{mf1} and \eqref{mf2} come from
estimating the solution $S(t)u_{0}$ of
the linear 1D Schr\"odinger IVP in the norm appearing in the
left hand side and a standard argument of summation along parabolic
curves, see for example the expository paper \cite{g}.
The smoothing effect and maximal function estimates
for $S(t)u_{0}$ can be found for example in \cite{kpv:s}.
\end{proof}

We also have the following improved Strichartz estimate
(cf. Lemma 7.1 in \cite{ckstt:0},\cite{borg:refinements} and
\cite{ot}):
\begin{lemma}\label{improved-strichartz}
For any Schwartz functions $u, v$ with Fourier support in $|\xi| \sim
R$, $|\xi| \ll R$ respectively, we have that
$$ \| u v \|_{L^2_t L^2_x} = \| u \bar v \|_{L^2_{t} L^2_x}
\lesssim R^{-1/2} \|u\|_{0,1/2+} \|v\|_{0,1/2+}.$$
\end{lemma}

In our arguments we shall be using the trivial embedding
$$ \| u \|_{s_1,b_1} \lesssim \| u \|_{s_2,b_2} \hbox{ whenever } s_1
\leq s_2, b_1 \leq b_2$$
so frequently that we will not mention this embedding explicitly.

We now give some useful notation for multilinear expressions.  If $n
\geq 2$ is an even integer, we define a \emph{(spatial)
$n$-multiplier} to be any function $M_n(\xi_1, \ldots, \xi_n)$ on the
hyperplane
$$ \Gamma_n := \{ (\xi_1, \ldots, \xi_n) \in \R^n: \xi_1 + \ldots +
\xi_n = 0 \},$$
which we endow with the standard measure $\delta(\xi_1 + \ldots +
\xi_n)$, where $\delta$ is the Dirac delta.

If $M_n$ is a $n$-multiplier and $f_1, \ldots, f_n$ are
functions on $\dbR$, we define the \emph{$n$-linear functional}
$\Lambda_n(M_n;f_1, \ldots,
f_n)$ by
$$ \Lambda_n(M_n;f_1, \ldots, f_n)
:= \int_{\Gamma_n} M_n(\xi_1, \ldots, \xi_n) \prod_{j=1}^n \hat
f_j(\xi_j).$$
We adopt the notation
$$ \Lambda_n(M_n;f) := \Lambda_n(M_n; f, \bar{f}, f, \bar{f}, \ldots,
f, \bar{f}).$$
Observe that $\Lambda_n(M_n;f)$ is invariant under permutations of
the even $\xi_j$ indices, or of the odd $\xi_j$ indices.

If $M_n$ is a multiplier of order $n$, $1 \leq j \leq n$ is an index,
and $k \geq 1$ is an even integer, we define the {\em elongation}
$\X^k_j(M_n)$ of $M_n$ to be the multiplier of order $n+k$ given by
$$ \X^k_j(M_n)(\xi_1, \ldots, \xi_{n+k})
:= M_n(\xi_1, \ldots, \xi_{j-1}, \xi_j + \ldots + \xi_{j+k},
\xi_{j+k+1}, \ldots, \xi_{n+k}).$$
In other words, $\X^k_j$ is the multiplier obtained by replacing
$\xi_j$ by $\xi_j + \ldots + \xi_{j+k}$ and advancing all the indices
after $\xi_j$ accordingly.

We shall often write $\xi_{ij}$ for $\xi_i + \xi_j$, $\xi_{ijk}$ for
$\xi_i + \xi_j + \xi_k$, etc.  We also write $\xi_{i-j}$ for $\xi_i -
\xi_j$, $\xi_{ij-klm}$ for $\xi_{ij} - \xi_{klm}$, etc. Also if
$m(\xi)$ is a function defined in the frequency space, we use the
notation $m(\xi_{i})=m_{i}, \, m(\xi_{ij-k})=m_{ij-k}$, etc.

In this paper we often use two very elementary tools: The mean value
theorem (MVT) and the double mean value theorem
(DMVT). While recalling
the statement of the MVT will be  an embarrassment, we think
that doing so for the DMVT  is a necessity to avoid later confusion.
\begin{lemma}[DMVT]
  Assume $f \in C^{2}(\dbR)$ and that $\max(|\eta|,|\lambda|)
  \ll |\xi|$, then
  $$|f(\xi+\eta+\lambda)-f(\xi+\eta)-f(\xi+\lambda)+f(\xi)|
  \lesssim |f''(\theta)||\eta||\lambda|,$$
  where $|\theta|\sim |\xi|$.
\label{dmvt}\end{lemma}

\section{The Gauge Transformation, energy and the almost
conservation  laws}
In this section we summarize the main results presented in Section 3
and 4 of \cite{ckstt:0}. Whatever is here simply  stated and
recalled is fully explained or proved in those sections.

We start by  applying the gauge transform
used in \cite{ozawa} in order to improve the derivative nonlinearity
present in \eqref{ivp}.

\begin{definition}\label{gauge-def}  We define the non-linear map
${\cal G}: L^2(\R) \to L^2(\R)$ by
$$
{\cal G} f(x) := e^{-i\int_{-\infty}^{x}|f(y)|^{2}dy}f(x).$$
The inverse transform ${\cal G}^{-1} f$ is then given by
$$
{\cal G}^{-1} f(x) := e^{i\int_{-\infty}^{x}|f(y)|^{2}dy}f(x).$$
\end{definition}

This transform is a bicontinuous map from $H^s$ to itself for any
$s\in [0,1]$.

Set $w_0 := {\cal G} u_0$, and $w(t) := {\cal G} u(t)$ for all times
$t$.  A straightforward calculation shows that the IVP \eqref{ivp}
  transforms into
\begin{equation}
\left\{ \begin{array}{l}
i\partial_tw+\partial_{x}^{2}w=-i w^{2}\partial_{x}\bar{w}
-\frac{1}{2}|w|^{4}w,\\
     w(x,0) = w_{0}(x),\hspace{1.5cm}x \in \dbR, \, t \in \dbR.
\end{array}\right.
\label{givp1}\end{equation}
Also, the smallness condition \eqref{small-l2} becomes
\begin{equation}\label{small-l2-w}
\| w_0 \|_{L^2} < \sqrt{2\pi}.
\end{equation}
By the bicontinuity we thus see that global well-posedness of
\eqref{ivp} in $H^s$ is equivalent to that of \eqref{givp1}.  From
\cite{ozawa, takaoka:dnls-local, takaoka:dnls-global},  we know that
both Cauchy problems are
locally well-posed in $H^s, s\geq \frac12$ and globally well-posed in
$H^1$ assuming
\eqref{small-l2-w}.  By standard limiting arguments, we thus see that
Theorem \ref{main} will follow if we can show:

\begin{proposition}\label{global}  Let $w$ be a global $H^1$ solution
to \eqref{givp1} obeying \eqref{small-l2-w}.  Then for any $T > 0$
and $s>\frac12$ we have
$$ \sup_{0 \leq t \leq T} \| w(t) \|_{H^s} \lesssim
C_{(\|w_0\|_{H^s},T)}$$
where the right-hand side does not depend on the $H^1$ norm of $w$.
\end{proposition}
We now pass to the considerations on the energy associated to
solutions of \eqref{givp1}.
\begin{definition}\label{energy-def}  If $f \in H^1(\R$), we define
the energy $E(f)$ by
$$ E(f) := \int \partial_x f\partial_{x}\overline{f} \ dx
- \frac{1}{2}\Im
\int f \overline{f} f \partial_x \overline{f}\ dx.$$
\end{definition}
By the Gagliardo-Nirenberg inequality we have
\begin{equation}\label{gee}
\|\partial_{x}f\|_2\leq C_{\|f\|_2} E(f)^{1/2},
\end{equation}
for any $f \in H^{1}$ such that $\|f\|_{2}<\sqrt{2\pi}$.

By Plancherel, we write $E(f)$ using the $\Lambda$ notation
and Fourier transform properties as
\begin{equation}
     E(f) = -\Lambda_2(\xi_1\xi_2;f) -  \frac{1}{2} \Im
\Lambda_4(i\xi_4;f).
\label{lambdaen}\end{equation}
Expanding out the second term using $\Im(z) = (z - \bar{z})/2i$, and
using symmetry, we may rewrite this as
\begin{equation}\label{energy-lambda}
E(f) = -\Lambda_2(\xi_1\xi_2;f) + \frac{1}{8} \Lambda_4(\xi_{13-24};f).
\end{equation}
One can use the same notation to rewrite the $L^{2}$ norm as
$$ \| w(t)\|_2^2 = \Lambda_2(1; w(t)).$$
\begin{lemma}\label{conservation}\cite{ozawa}  If $w$ is an $H^1$
solution to \eqref{givp1} for $t \in [0,T]$, then we have
$$ \| w(t) \|_2 = \| w_0 \|_2$$
and
$$ E(w(t)) = E(w_0)$$
for all $t \in [0,T]$.
\end{lemma}
In \cite{ckstt:0} this lemma was proved using the following general
proposition (cf. \cite{ckstt:0}):
\begin{proposition}
     Let $n\geq 2$ be an even integer, let $M_{n}$ be a multiplier of
     order $n$ and let $w$ be a solution of \eqref{givp1}. Then
     \begin{equation}\label{diff}
\begin{split}
\partial_t \Lambda_n(M_n;w(t))
&=
i \Lambda_n(M_n \sum_{j=1}^n (-1)^{j} \xi_j^2; w(t))\\
&- i \Lambda_{n+2}(\sum_{j=1}^n \X^2_j(M_n) \xi_{j+1}; w(t))\\
&+\frac{i}{2} \Lambda_{n+4}(\sum_{j=1}^n (-1)^{j-1} \X^4_j(M_n);
w(t)).
\end{split}
\end{equation}
\label{ddt}\end{proposition}
We summarize below the idea we used to prove Proposition
\ref{global}  for $s>\frac{2}{3}$ in \cite{ckstt:0}.
Because we do not want to use the $H^1$ norm of $w$, we cannot
directly use the energy $E(w(t))$ defined
above. So we introduced a substitute notion of ``energy''
that could be
defined for a less regular solution and that had a very slow increment
in time. In frequency space consider an even  $C^{\infty}$
monotone multiplier
$m(\xi)$ taking values in $[0,1]$ such that
\begin{equation}
m(\xi):=\left\{ \begin{array}{l}
1, \, \, \, \mbox{ if } \, |\xi|<N,\\
\left(\frac{|\xi|}{N}\right)^{s-1}\, \, \, \mbox{ if } \, |\xi|> 2N.
\end{array}\right.
\label{mul}\end{equation}
Define the multiplier operator $I:H^{s}\longrightarrow H^{1}$
such that $\widehat{Iw}(\xi):=m(\xi)\widehat{w}(\xi)$.  This operator
is smoothing of order $1-s$; indeed one has
\begin{equation}\label{i-smoothing}
\| u\|_{s_0,b_0} \lesssim \| Iu \|_{s_0 + 1-s, b_0} \lesssim N^{1-s}
\| u \|_{s_0,b_0}
\end{equation}
for any $s_0, b_0 \in \R$. Our substitute energy was defined by
$$ E_N(w) := E(Iw).$$
Note that this energy makes sense even if $w$ is only in $H^s$.
In general the energy $E_N(w(t))$ is not conserved in time, but we
showed that the increment was very small in terms of $N$.

To proceed with the improvement of the ``I-method'',
let us consider a
symmetric multiplier $m(\xi)$\footnote{This eventually will be taken
to be exactly the multiplier in \eqref{mul}.} and let $I$ be the
multiplier operator associated to it. Then we write
$$E^{1}(w):= E(Iw).$$
Clearly, if $m$ is the multiplier in \eqref{mul}, then
$$E^{1}(w)=E_N(w), $$
so we can think  about $E^{1}(w)$  as the first generation of a
family of modified energies. In this paper we introduce the second
generation in detail, but formally the method can be used to
define an infinite family of modified energies. We write
\begin{equation}
     E^{2}(w)=-\Lambda_2(m_1\xi_1m_2\xi_2,w)
+\frac12\Lambda_4\left(M_4(\xi_1,\xi_2,\xi_3,\xi_4),w\right),
     \label{gen2}\end{equation}
where $M_{4}$ will be determined  later. Assume now that $w$ is
a solution
of \eqref{givp1}. Because  $w$ is fixed we drop it from the
definition of $E^{2}$. We are interested in the
increment of this second generation of energies, hence we compute
$\frac{d}{dt}E^{2}$.
Differentiating $\Lambda_2(m_1\xi_1m_2\xi_2)$
using Proposition \ref{ddt}, using the identity $\xi_{1}+\ldots+
\xi_{n}=0$ and symmetrizing, we have
\begin{eqnarray*}
\frac{d}{dt}\Lambda_2(m_1\xi_1m_2\xi_2)
&=&
-i\Lambda_2(m_1\xi_1m_2\xi_2(\xi_1^2-\xi_2^2))
-i\Lambda_4(m_{123}\xi_{123}m_4\xi_4\xi_2+m_1\xi_1m_{234}
\xi_{234}\xi_3)\\
& &
+\frac{i}{2}\Lambda_6(m_{12345}\xi_{12345}m_6\xi_6-m_1\xi_1m_{23456}
\xi_{23456})\\
&=&
\frac{i}{2}\Lambda_4(\sigma_4(\xi_1,\xi_2,\xi_3,\xi_4))+
\frac{i}{6}\Lambda_6(\sigma_6(\xi_1,\xi_2,\xi_3,\xi_4,\xi_5,\xi_6)),
\end{eqnarray*}
where
\begin{equation}
\sigma_4(\xi_1,\xi_2,\xi_3,\xi_4)=m_1^2\xi_1^2\xi_3+m_2^2\xi_2^2
\xi_4+m_3^2\xi_3^2\xi_1+m_4^2\xi_4^2\xi_2,
\label{sigma4}\end{equation}
and
\begin{equation}
\sigma_6(\xi_1,\xi_2,\xi_3,\xi_4,\xi_5,\xi_6)=\sum_{j=1}^6(-1)^{j-1}
m_j^2\xi_j^2.
\label{sigma6}\end{equation}
Notice that the contribution of $\Lambda_{2}$ is zero
because the factor $(\xi_1^2-\xi_2^2)$ is zero over the set of
integration $\xi_1+\xi_2=0$.

Differentiating $\Lambda_4(M_4)$, we have
\begin{eqnarray*}
& &
\frac{d}{dt}\Lambda_4(M_4(\xi_1,\xi_2,\xi_3,\xi_4))\\
&=&
i\Lambda_4(M_4\sum_{j=1}^{4}(-1)^{j}\xi_{j}^{2})\\
& &
-i\Lambda_6(M_4(\xi_{123},\xi_4,\xi_5,\xi_6)\xi_2+M_4(\xi_1,
\xi_{234},\xi_5,\xi_6)\xi_3+M_4(\xi_1,\xi_2,\xi_{345},\xi_6)
\xi_4+M_4(\xi_1,\xi_2,\xi_3,\xi_{456})\xi_5)\\
& &
+\frac{i}{2}\Lambda_8(M_4(\xi_{12345},\xi_6,\xi_7,\xi_8)-
M_4(\xi_1,\xi_{23456},\xi_7,\xi_8)+M_4(\xi_1,\xi_2,\xi_{34567},
\xi_8)-M_4(\xi_1,\xi_2,\xi_3,\xi_{45678}))\\
&=&
i\Lambda_4(M_4\sum_{j=1}^{4}(-1)^{j}\xi_{j}^{2})\\
& &
-\frac{i}{36}\sum_{{\scriptstyle\{a,c,e\}=\{1,3,5\}}
\atop{\scriptstyle\{b,d,f\}=\{2,4,6\}}}\Lambda_6(M_4(\xi_{abc},
\xi_d,\xi_e,\xi_f)\xi_b+M_4(\xi_a,\xi_{bcd},\xi_e,\xi_f)\xi_c\\
&
&+M_4(\xi_a,\xi_b,\xi_{cde},\xi_f)\xi_d+M_4(\xi_a,\xi_b,\xi_c,
\xi_{def})\xi_e)\\
& &
+C\sum_{{\scriptstyle\{a,c,e,g\}=\{1,3,5,7\}}
\atop{\scriptstyle\{b,d,f,h\}=\{2,4,6,8\}}}\Lambda_8(M_4(\xi_{abcde},
\xi_f,\xi_g,\xi_h)+M_4(\xi_a,\xi_b,\xi_{cdefg},\xi_h)\\
& &
-M_4(\xi_a,\xi_{bcdef},\xi_g,\xi_h)-M_4(\xi_a,\xi_b,\xi_c,
\xi_{defgh}))
\end{eqnarray*}
Then
\begin{eqnarray*}
\frac{d}{dt}E^{2}(w)&=&
-\frac{i}{2}\Lambda_4(\sigma_{4}(\xi_{1},\xi_{2},\xi_{3},\xi_{4}))
+\frac{i}{2}\Lambda_4(M_4\sum_{j=1}^{4}(-1)^{j}\xi_{j}^{2})\\
&-&\frac{i}{6}
\Lambda_6(\sigma_6(\xi_1,\xi_2,\xi_3,\xi_4,\xi_5,\xi_6))\\
& &-\frac{i}{72}\sum_{{\scriptstyle\{a,c,e\}=\{1,3,5\}}
\atop{\scriptstyle\{b,d,f\}=\{2,4,6\}}}\Lambda_{6}(M_4(\xi_{abc},
\xi_d,\xi_e,\xi_f)\xi_b+M_4(\xi_a,\xi_{bcd},\xi_e,\xi_f)
\xi_c\\
& &+M_4(\xi_a,\xi_b,\xi_{cde},\xi_f)\xi_d+M_4(\xi_a,
\xi_b,\xi_c,\xi_{def})\xi_e)\\
& &
+C_1\sum_{{\scriptstyle\{a,c,e,g\}=\{1,3,5,7\}}
\atop{\scriptstyle\{b,d,f,h\}=\{2,4,6,8\}}}
\Lambda_8(M_4(\xi_{abcde},\xi_f,\xi_g,\xi_h)+
M_4(\xi_a,\xi_b,\xi_{cdefg},\xi_h)\\
& &
-M_4(\xi_a,\xi_{bcdef},\xi_g,\xi_h)-M_4(\xi_a,\xi_b,
\xi_c,\xi_{defgh})).
\end{eqnarray*}
We abbreviate the 6-linear and the 8-linear expressions as
$\Lambda_6(M_6(\xi_1,\xi_2,\cdots,\xi_6))$ and
$\Lambda_8(M_8(\xi_1,\xi_2,\cdots,\xi_8))$.
We are now ready to make our choice for $M_{4}$. From our
calculations in \cite{ckstt:0}, we realized that the estimates for
the different pieces of $\Lambda_{n}$ appearing in the right had side
if $\frac{d}{dt}E_{N}(w)$ are easier  for $n$ larger\footnote{Compare
for example section 8, 9 and 10 in \cite{ckstt:0}.}, we decided to
use the freedom of choosing $M_{4}$ to cancel the $\Lambda_{4}$
contribution obtained above. Hence using \eqref{sigma4}, we set
\begin{equation}
M_4(\xi_1,\xi_2,\xi_3,\xi_4)=-\frac{m_1^2\xi_1^2\xi_3+m_2^2
\xi_2^2\xi_4+m_3^2\xi_3^2\xi_1+m_4^2\xi_4^2\xi_2}
{\xi_1^2-\xi_2^2+\xi_3^2-\xi_4^2},
\label{m4}\end{equation}
which in the set of integration $\xi_{1}+\xi_{2}+\xi_{3}+\xi_{4}=0$,
can also be written as
$$M_4(\xi_1,\xi_2,\xi_3,\xi_4)=
-\frac{m_1^2\xi_1^2\xi_3+m_2^2\xi_2^2\xi_4+m_3^2\xi_3^2
\xi_1+m_4^2\xi_4^2\xi_2}{2\xi_{12}\xi_{14}}.$$
\begin{remark}
  If we assume that $m(\xi)=1$, then $E^{2}(w)=E(w)$. In fact, on the
  set $\xi_{1}+\xi_{2}+\xi_{3}+\xi_{4}=0$  we have
\begin{eqnarray*}
& &m_1^2\xi_1^2\xi_3+m_2^2
\xi_2^2\xi_4+m_3^2\xi_3^2\xi_1+m_4^2\xi_4^2\xi_2=
\xi_1^2\xi_3+
\xi_2^2\xi_4+\xi_3^2\xi_1+\xi_4^2\xi_2
=(\xi_{1}+\xi_{3})(\xi_{1}\xi_{3}-\xi_{2}\xi_{4})\\
&=&(\xi_{1}+\xi_{3})(\xi_{1}\xi_{3}+(\xi_{1}+\xi_{3}+\xi_{4})\xi_{4})
=-(\xi_{1}+\xi_{3})(\xi_{1}+\xi_{4})(\xi_{1}+\xi_{2}),
\end{eqnarray*}
hence
\begin{equation}
\label{m41}
M_4(\xi_1,\xi_2,\xi_3,\xi_4)=\frac{1}{2}(\xi_{1}+\xi_{3}),
\end{equation}
and
$$E^{2}(w)=-\Lambda_2(\xi_1\xi_2)+\frac14\Lambda_4(\xi_{13})$$
which is exactly the value of $E(w)$ in \eqref{lambdaen}.
\label{remm4}\end{remark}
Once again we recall that we assume throughout the paper
that $s\in (\frac12, \frac23]$ and that the
multiplier $m$ is defined as in \eqref{mul}. To stress the fact
that with this choice the energy $E^{2}(w)$  depends on the
parameter $N$, we  write $E^{2}(w)=E^{2}_{N}$. We now summarize
some of the  above observations in the following:
\begin{proposition}\label{energy-increment}
Let $w$ be an $H^1$ global solution to \eqref{givp1}.  Then for any
$T \in \R$ and $\delta > 0$ we have
$$
E^{2}_N(w(T+\delta)) - E^{2}_N(w(T)) =
\int_T^{T+\delta} [\Lambda_6(M_6;w(t)) +\Lambda_8(M_8;w(t))]\ dt
$$
where the multipliers $M_6$ and $M_8$ are given by
\begin{eqnarray*}
M_6 &:=& -\frac{i}{6}\sigma_6(\xi_1,\xi_2,
\xi_3,\xi_4,\xi_5,\xi_6)\\
& &-\frac{i}{72}\sum_{{\scriptstyle\{a,c,e\}=\{1,3,5\}}
\atop{\scriptstyle\{b,d,f\}=\{2,4,6\}}}(M_4(\xi_{abc},
\xi_d,\xi_e,\xi_f)\xi_b+M_4(\xi_a,\xi_{bcd},\xi_e,\xi_f)
\xi_c\\
& &+M_4(\xi_a,\xi_b,\xi_{cde},\xi_f)\xi_d+M_4(\xi_a,
\xi_b,\xi_c,\xi_{def})\xi_e) \\
M_8 &:=& C_2\sum_{{\scriptstyle\{a,c,e,g\}=\{1,3,5,7\}}
\atop{\scriptstyle\{b,d,f,h\}=\{2,4,6,8\}}}
(M_4(\xi_{abcde},\xi_f,\xi_g,\xi_h)+
M_4(\xi_a,\xi_b,\xi_{cdefg},\xi_h)\\
& &
-M_4(\xi_a,\xi_{bcdef},\xi_g,\xi_h)-M_4(\xi_a,\xi_b,
\xi_c,\xi_{defgh}))
\end{eqnarray*}
where $C_2$ is an absolute constant.
Furthermore, if $|\xi_j| \ll N$ for all $j$, then the multipliers
$M_6$ and $M_8$  vanish.
\end{proposition}
We end this section with a lemma that shows
the energy $E^{2}_{N}(w)$ has the same strength as $\|Iw\|_{H^{1}}$.
\begin{lemma}
Assume that $w$  satisfies
$\|w\|_{L^2}<\sqrt{2\pi},~\|Iw\|_{H^1}=O(1)$.
Then for $N>>1$,
\begin{equation}
\|\partial_xIw\|_{L^2}^{2}\lesssim E^{2}_{N}(w).
\label{E2I}\end{equation}
\label{dxe}\end{lemma}
The proof of this lemma relies strongly on the estimate of the
multiplier $M_{4}$ and it can be found in the next section.

\section{Estimates for $M_{4}$ and proof of Lemma \ref{dxe}}
Before we start with our estimates
we recall some notation that we used in \cite{ckstt:0}.
Let $n = 4$, $6$, or $8$, and let $\xi_1,
\ldots, \xi_n$ be frequencies such that $\xi_1 + \ldots + \xi_n =
0$.  Define $N_i := |\xi_i|$, and $N_{ij} := |\xi_{ij}|$.  We adopt
the notation that
$$ 1 \leq soprano, alto, tenor, baritone \leq n$$
are the distinct indices such that
$$ N_{soprano} \geq N_{alto} \geq N_{tenor} \geq N_{baritone}$$
are the highest, second highest, third highest, and fourth highest
values of the frequencies $N_1, \ldots, N_n$ respectively (if there
is a tie in frequencies, we break the tie arbitrarily).
Since $\xi_1 + \ldots + \xi_n = 0$, we must have $ N_{soprano} \sim
N_{alto}$.  Also, from Proposition \ref{energy-increment} we see that
$M_n$ vanishes unless $N_{soprano} \gtrsim N$.

In this section whenever we write $\max|f(\theta)|$,
for a function $f$ we understand that the maximum is taken for
$|\theta|\sim N_{soprano}$.

\begin{lemma}
Assume $M_{4}$ is the multiplier defined in \eqref{m4} and $m(\xi)$ is
like in \eqref{mul}. Then
\begin{equation}
     |M_{4}(\xi_{1},\ldots,\xi_{4})|\lesssim
m^{2}(N_{soprano})N_{soprano}.
\label{m4est}\end{equation}
\label{estm4}
\end{lemma}
\begin{proof}
We  observe that to prove \eqref{m4est} it suffices to prove
$$|\sigma_4(\xi_{1},\ldots,\xi_{4})|\lesssim |\xi_{12}||\xi_{12}|
m^{2}(N_{soprano})N_{soprano}.$$
Without loss of generality we may assume that
$N_{soprano}=N_{1}$. By symmetry we can assume that
$|\xi_{12}|\le|\xi_{14}|$.
We divide the analysis into two cases: Case a) when
$N_{1}\lesssim |\xi_{14}|$ and Case b) when $|\xi_{14}|\ll N_{1}$.

\noindent
{\bf Case a):} we write
\begin{eqnarray}
\nonumber|\sigma_4(\xi_{1},\ldots,\xi_{4})|&=&
|m_1^2\xi_1^2\xi_3+m_2^2\xi_2^2(-\xi_{12}-\xi_3)+
m_3^2\xi_3^2\xi_1+m_{12+3}^2\xi_{12+3}^2\xi_2|\\
  \label{rewrite}&=&
|\xi_3(m_1^2\xi_1^2-m_{1-12}^2\xi_{1-12}^2)+\xi_1(m_3^2\xi_3-
m_{3+12}^2\xi_{3+12}^2)-\xi_{12}(m_2^2\xi_2^2-m_{12+3}^2
\xi_{12+3}^2)|.
\end{eqnarray}
Then the MVT shows that
\begin{equation}
|\sigma_4(\xi_1,\xi_2,\xi_3,\xi_4)|\lesssim
|\xi_{12}|N_{1}\max|(m(\xi)^{2}\xi^2)'|,
\label{eq:MV}\end{equation}
where $|\xi|\lesssim N_{1}$.
Now it is easy to see that for $m$ defined in \eqref{mul}
$$(m^{2}(\xi)\xi^{2})'\sim m^{2}(\xi)\xi,$$
and that the function $m^{2}(\xi)\xi$ is non decreasing.
Then  \eqref{eq:MV} immediately gives \eqref{m4est}.

\noindent
{\bf Case b):} We first write $\sigma_{4}$ so that
the DMVT in Lemma \ref{dmvt} can be applied. For simplicity we write
$m^{2}(\xi)\xi^{2}=f(\xi)$. Then in the set
$\xi_{1}+\ldots+ \xi_{4}=0$ we have
\begin{eqnarray*}
\sigma_{4}(\xi_{1},\ldots, \xi_{4})&=&f(\xi_{1})\xi_{3}
+f(\xi_2)\xi_{4}+f(\xi_{3})\xi_{1}f(\xi_{4})\xi_{2}\\
&=&
\xi_3[f(\xi_{1})-f(\xi_2)]+\xi_1[f(\xi_3)-
f(-\xi_{4})]-\xi_{12}[f(\xi_2)-f(-\xi_{4})]\\
&=&\xi_{3}[f(\xi_{1})-f(\xi_{2})+f(\xi_{3})-f(-\xi_{4})]
\\&+&(\xi_{1}-\xi_{3})
[f(\xi_{3})-f(\xi_{3}-\xi_{12})]-\xi_{12}[f(\xi_{2})-f(-\xi_{4})]\\
&=&\xi_{3}[f(\xi_{1}-\xi_{12}-\xi_{14})-f(\xi_{1}-\xi_{12})
-f(\xi_{1}-\xi_{14})+f(\xi_{1})]\\
&+&(-\xi_{3}+\xi_{1})
[f(\xi_{3})-f(\xi_{3}-\xi_{12})]-\xi_{12}[f(\xi_{2})-f(\xi_{2}
+\xi_{14}-\xi_{12})].
\end{eqnarray*}
where we often used the fact that $f(\xi)$ is an even function.
Using the DMVT in the first term of the right hand side of the
inequality and the MVT in the remaining two terms we obtain
\begin{equation}
\sigma_{4}(\xi_{1},\ldots, \xi_{4})\lesssim
|\xi_{1}||f''(\theta)||\xi_{12}||\xi_{14}| +
|\xi_{12}|\max|f'|(|\xi_{3-1}|+|\xi_{14}|+|\xi_{12}|),
     \label{step1}\end{equation}
where $|\theta|\sim N_{1}$. Now observe that
$$|\xi_{3-1}|=|\xi_{12}+\xi_{14}|\lesssim |\xi_{14}|$$
and that  $|f''(\theta)|\lesssim m(N_{1})^{2}$,
so inserting \eqref{step1} in the definition of $M_{4}$ we obtain
\eqref{m4est}.
\end{proof}
We need two more local estimates for $M_{4}$:
\begin{lemma}

\noindent
\begin{itemize}
     \item Assume that $|\xi_1|\sim|\xi_3|\gtrsim N\gg |\xi_2|,|\xi_4|$,
then
\begin{equation}
|M_4(\xi_1,\xi_2,\xi_3,\xi_4)|\lesssim m(N_{soprano})^2N_{tenor}.
\label{13m4}\end{equation}
\item
Assume that $|\xi_1|\sim|\xi_2|\gtrsim N\gg |\xi_3|,|\xi_4|$, then
\begin{equation}
M_4(\xi_1,\xi_2,\xi_3,\xi_4)=\frac{m_1^2\xi_2^2}{2\xi_1}+
R(\xi_{1},\ldots,\xi_{4}),
\label{m4app}\end{equation}
where
$$|R(\xi_{1},\ldots,\xi_{4})|\lesssim N_{tenor}.$$
\end{itemize}
\end{lemma}
\begin{proof}
The  first part of the lemma  follows from the MVT. In fact
\begin{eqnarray*}
\left|\frac{m_1^2\xi_1^2\xi_3+\xi_2^2\xi_4+m_3^2
\xi_3^2\xi_1+\xi_4^2\xi_2}{\xi_{12}\xi_{14}}
\right|\lesssim\frac{|\xi_1\xi_3\xi_{13}|\max|(m(\xi)^{2}\xi)'|+
|\xi_{24}\xi_2\xi_4|}{|\xi_1|^2}\lesssim
m(N_{soprano})^2N_{tenor},
\end{eqnarray*}
where again we used that $|(m(\xi)^{2}\xi)'|\sim |m(\xi)\xi|$.

To prove the second part of the lemma we use the identity
$$\frac{1}{\xi_{14}}=\frac{1}{\xi_{1}}-\frac{\xi_{4}}{\xi_{14}}
\frac{1}{\xi_{1}},$$
and we write
$$-2M_4(\xi_1,\xi_2,\xi_3,\xi_4)+\frac{m_1^2\xi_2^2}{\xi_1}=
R_{1}(\xi_{1},\ldots,\xi_{4})+R_{2}(\xi_{1},\ldots,\xi_{4}),$$
where
\begin{eqnarray*}
  R_{1}(\xi_{1},\ldots,\xi_{4})&=&\frac{m_{1}^{2}\xi_{1}^{2}\xi_{3}
     +m_{2}^{2}\xi_{2}^{2}\xi_{4}+\xi_{3}^{2}\xi_{1}
     +\xi_{4}^{2}\xi_{2}+m_{1}^{2}\xi_{2}^{2}\xi_{12}}
     {\xi_{12}\xi_{1}}\\
     R_{2}(\xi_{1},\ldots,\xi_{4})&=&-\frac{\xi_{4}}{\xi_{14}}
     \frac{m_{1}^{2}\xi_{1}^{2}\xi_{3}
     +m_{2}^{2}\xi_{2}^{2}\xi_{4}+\xi_{3}^{2}\xi_{1}
     +\xi_{4}^{2}\xi_{2}}{\xi_{12}\xi_{1}}.
\end{eqnarray*}
We estimate first $R_{1}$:
\begin{eqnarray*}
R_{1}(\xi_{1},\ldots,\xi_{4})&=&\frac{m_{1}^{2}\xi_{1}^{2}\xi_{3}
     +m_{2}^{2}\xi_{2}^{2}\xi_{4}+\xi_{3}^{2}\xi_{1}
     +\xi_{4}^{2}\xi_{2}-m_{1}^{2}\xi_{2}^{2}\xi_{34}}
     {\xi_{12}\xi_{1}}\\
  &=&\frac{m_{1}^{2}\xi_{3}(\xi_{1}^{2}-\xi_{2}^{2})+
  \xi_{2}^{2}\xi_{4}(m_{2}^{2}-m_{1}^{2})+
  \xi_{3}^{2}(\xi_{1}+\xi_{2})+\xi_{2}(\xi_{4}^{2}-\xi_{3}^{2})}
  {\xi_{12}\xi_{1}},
  \end{eqnarray*}
hence, by the MVT,
$$|R_{1}(\xi_{1},\ldots,\xi_{4})|\lesssim N_{tenor}.$$
On the other hand
\begin{eqnarray*}
R_{2}(\xi_{1},\ldots,\xi_{4})&=&-\frac{\xi_{4}}{\xi_{14}}
\frac{m_{1}^{2}\xi_{1}^{2}(\xi_{3}+\xi_{4})+
     (m_{2}^{2}\xi_{2}^{2}-m_{1}^{2}\xi_{1}^{2})\xi_{4}+
  \xi_{3}^{2}\xi_{12}+\xi_{2}\xi_{34}\xi_{3-4}}
  {\xi_{12}\xi_{1}},
  \end{eqnarray*}
hence, again by the MVT,
$$|R_{2}(\xi_{1},\ldots,\xi_{4})|\lesssim N_{tenor}.$$
\end{proof}

\noindent
{\bf Proof of Lemma \ref{dxe}}
\begin{proof}
We rewrite $E^{2}_{N}(w)$ as
\begin{eqnarray*}
E^{2}_{N}(w)&=&-\Lambda_2(m_1\xi_1m_2\xi_2)+\frac18\Lambda_4
(\xi_{13-24}m_1m_2m_3m_4)\\
& &+\frac18\Lambda_4(4M_4(\xi_1,\xi_2,\xi_3,\xi_4)
-\xi_{13-24}m_1m_2m_3m_4).
\end{eqnarray*}
In Lemma 3.6 of \cite{ckstt:0} we proved the estimate
$$
\|\partial_xIw\|_{L^2}^{2}\lesssim -\Lambda_2(m_1\xi_1m_2\xi_2)+
\frac18\Lambda_4(\xi_{13-24}m_1m_2m_3m_4)
$$
for $\|Iw\|_{L^2}<\sqrt{2\pi}$.
Hence we only have  to show that
\begin{equation}\label{eq:Ga}
\left|\Lambda_4(4M_4(\xi_1,\xi_2,\xi_3,\xi_4)-
\xi_{13-24}m_1m_2m_3m_4)\right|\lesssim
O\left(\frac{1}{N^{\alpha}}\right)\|Iw\|_{H^1}^4
\end{equation}
for some $\alpha>0$.

We first perform a Littlewood-Paley decomposition of the four
factors $w$ so that the $\xi_{i}$ are essentially the constants
$N_{i}, \, i=1,\ldots,4$.
To recover the sum at the end we borrow a $N^{-\epsilon}_{soprano}$
from the
large denominator $N_{soprano}$ and often this will not be mentioned.

If all $|\xi_j|$ are less than $\frac{N}{100}$, the left hand side of
(\ref{eq:Ga}) vanishes thanks to \eqref{m4}.
Therefore, we may assume  $N_{soprano}\gtrsim N$.
Also note $N_{alto}\gtrsim N$ on the set
$\xi_{1}+\xi_{2}+\xi_{3}+\xi_{4}=0$. Then it is obvious that
$$
|\Lambda_4(\xi_{13-24}m_1m_2m_3m_4)|\lesssim\frac{1}{N}
\|Iw\|_{H^1}^2\|Iw\|_{L^{\infty}}^2\lesssim\frac{1}{N}\|Iw\|_{H^1}^4.
$$
Next we control the contribution of $\Lambda_4(M_4)$ in (\ref{eq:Ga}).
By (\ref{m4est}), we have
\begin{eqnarray*}
|\Lambda_4(M_4(\xi_1,\xi_2,\xi_3,\xi_4))|
\lesssim\frac{1}{N_{soprano}^{1-}m(N_{baritone})^{2}
N_{baritone}}
\|Iw\|_{H^1}^4\lesssim\frac{1}{N^{1-}}\|Iw\|_{H^1}^4,
\end{eqnarray*}
where again we used the fact that $m^{2}(\xi)\xi$ is non decreasing.
\end{proof}

\section{Local Estimates}
This section contains a refinement of the results presented in
Section  5 of \cite{ckstt:0}. We start with the main result:
\begin{theorem}
Let $w$ be a $H^1$ global solution to \eqref{givp1} and let
$T \in \R$ be such that
$$ \| Iw(T) \|_{H^1} \leq C_0$$
for some $C_0 > 0$.
Then we have
\begin{eqnarray*}
\|Iw\|_{X^{1,b}([T,T+\delta]\times\Bbb R)}\lesssim 1
\end{eqnarray*}
for any $\frac12<b<\frac34$ and for some $\delta > 0$
depending on $C_0$.
\label{lwp}\end{theorem}
\begin{remark}
This theorem is stronger than the corresponding Theorem 5.1 in
\cite{ckstt:0} because $b$ can be arbitrarily close to $\frac34$,
and this is essential to obtain our sharp global well-posedness
result.
\end{remark}
As explained in \cite{ckstt:0} the proof of Theorem \ref{lwp} is a
consequence of the following multilinear estimates.
\begin{lemma}\label{lem:3.2}
For the Schwartz function $w$ and $\frac12<b<\frac34,~b'<\frac34$,
we have
\begin{eqnarray}\label{eq:Nt}
\|I(w\partial_x\overline{w}w)\|_{1,b'-1}\lesssim
\|Iw\|_{1,\frac12+}^2\|Iw\|_{1,b},
\end{eqnarray}
\begin{equation}\label{eq:Nt1}
\|I(w\overline{w}w\overline{w}w)\|_{1,b'-1}\lesssim
\|Iw\|_{1,\frac12+}^5.
\end{equation}
\end{lemma}
\begin{proof}
The proof of \eqref{eq:Nt1} follows from the same arguments used
to prove
(17) in \cite{ckstt:0}, and we do not  present it here again.
The proof of \eqref{eq:Nt} on the other hand is more delicate than the
one  given in \cite{ckstt:0} for (16), so we decided to give all the
details. By standard duality
arguments  in $L^2$ and renormalization,
it is easy to see that (\ref{eq:Nt}) is equivalent to
\begin{eqnarray}\label{eq:DNt}
\int_*
\frac{m_4\langle\xi_4\rangle|\xi_2|\langle\tau_4+\xi_4^2
\rangle^{b'-1}}
{\sum_{i=1}^3\langle\tau_i+(-1)^i\xi_i^2\rangle^{b-\frac12-}
\prod_{j=1}^3m_j\langle\xi_j\rangle\langle\tau_j+(-1)^j\xi_j^2
\rangle^{\frac12+}}
\prod_{j=1}^4F_j(\tau_j,\xi_j)\lesssim
\prod_{j=1}^4\|F_j\|_{L^2},
\end{eqnarray}
where all functions $F_j$ are real-valued and non-negative. If
\begin{eqnarray}\label{eq:R}
\frac{m_4\langle\xi_4\rangle|\xi_2|}{\prod_{j=1}^3m_j
\langle\xi_j\rangle}\lesssim 1,
\end{eqnarray}
then the $L^2$ estimate (\ref{strichartz-0}) for $F_4$ and
the  Strichartz estimate \eqref{strichartz-p} with $p=6$ for
$F_1,F_2,F_3$,
automatically shows (\ref{eq:DNt}) for $b>\frac12,~b'\le 1$.
Then we may assume
\begin{eqnarray*}
\frac{m_4\langle\xi_4\rangle|\xi_2|}{\prod_{j=1}^3m_j\langle
\xi_j\rangle}\gg 1,
\end{eqnarray*}
which, one can easily check, can happen only when
\begin{eqnarray*}
|\xi_2|\gg 1,~|\xi_{12}|\gg 1,~|\xi_{14}|\gg 1.
\end{eqnarray*}
We recall (cf. \cite{borg:xsb} and \cite{ckstt:0})
the fundamental inequality
\begin{equation}
     |\xi_{12}\xi_{14}|\lesssim \max_{j=1,2,3,4}
     \{\langle\tau_j+(-1)^j\xi_j^2\rangle\}.
\label{fund}\end{equation}
Then we  proceed with a case by case analysis: Case a) if
$\max_{j=1,2,3}\{\langle\tau_4+\xi_4^2\rangle, \langle
\tau_j+(-1)^j\xi_j^2\rangle\}=\langle\tau_4+\xi_4^2\rangle$ and
Case~b) if
$\max_{j=1,2,3}\{\langle\tau_4+\xi_4^2\rangle, \langle
\tau_j+(-1)^j\xi_j^2\rangle\}=\langle
\tau_i+(-1)^j\xi_i^2\rangle$, for some $i=1,2,3$.
\begin{itemize}
\item {\bf Case a):} In this case we replace in the denominator
$\langle\tau_4+\xi_4^2\rangle^{1-b'}$ with $(\langle\xi_{12}\rangle
\langle\xi_{14}\rangle)^{1-b'}$. Then using the same argument that
in \cite{ckstt:0}
led us from (16) to (18), we can show that \eqref{eq:DNt} is
equivalent to
\begin{eqnarray}\label{eq:DNt1}
\int_*
\frac{\langle\xi_4\rangle^{s}\langle\xi_2\rangle^{1-s}}
{(\langle\xi_{12}\rangle
\langle\xi_{14}\rangle)^{1-b'}
\langle\xi_1\rangle^{s}\langle\xi_3\rangle^{s}
\prod_{j=1}^3\langle\tau_j+(-1)^j\xi_j^2
\rangle^{\frac12+}}
\prod_{j=1}^4F_j(\tau_j,\xi_j)\lesssim
\prod_{j=1}^4\|F_j\|_{L^2}.
\end{eqnarray}
To have an idea of the ``numerics'' involved while proceeding with the
proof,
the reader should keep in mind that the interesting case is when
$s=\frac12+$ and $1-b'= \frac14+$. Since
$\xi_{14}=-\xi_{32}$, by symmetry, we may assume
that $|\xi_1|\ge|\xi_3|$. Then using the fact that
$\xi_{4}=-\xi_{3}-\xi_{12}$, we can write
\begin{equation}
\frac{\langle\xi_4\rangle^{s}\langle\xi_2\rangle^{1-s}}
{(\langle\xi_{12}\rangle
\langle\xi_{14}\rangle)^{1-b'}
\langle\xi_1\rangle^{s}\langle\xi_3\rangle^{s}}=A_{1}+A_{2},
\label{mult}\end{equation}
where
\begin{eqnarray*}
A_{1}&\lesssim&\frac{\langle\xi_2\rangle^{1-s}}
{(\langle\xi_{12}\rangle
\langle\xi_{14}\rangle)^{1-b'}
\langle\xi_1\rangle^{s}}\\
A_{2}&\lesssim&\frac{\langle\xi_{12}\rangle^{s-1+b'}
\langle\xi_2\rangle^{1-s}}
{\langle\xi_{14}\rangle^{1-b'}
\langle\xi_1\rangle^{s}\langle\xi_3\rangle^{s}}.
  \end{eqnarray*}
We now write $\xi_{12}=-\xi_{14}-\xi_{3}+\xi_{1}$ and we write
$$A_{2}=A_{2}^{1}+A_{2}^{2}+A_{2}^{3}$$
where
\begin{eqnarray*}
\\
A_{2}^{1}&\lesssim&\frac{\langle\xi_2\rangle^{1-s}}
{\langle\xi_{14}\rangle^{2(1-b')-s}
\langle\xi_1\rangle^{s}\langle\xi_3\rangle^{s}}\\
A_{2}^{2}&\lesssim&\frac{\langle\xi_2\rangle^{1-s}}
{\langle\xi_{14}\rangle^{1-b'}
\langle\xi_3\rangle^{1-b'}\langle\xi_1\rangle^{s}}\\
A_{2}^{3}&\lesssim&\frac{\langle\xi_2\rangle^{1-s}}
{\langle\xi_{14}\rangle^{1-b'}
\langle\xi_1\rangle^{1-b'}\langle\xi_3\rangle^{s}}.
\end{eqnarray*}
It is now easy to see that for $1-b'\geq \frac{s}{2}$,
$$A_{1}, A_{2}^{i}(\xi_{1},\xi_{2},\xi_{3})\lesssim
\frac{\langle\xi_2\rangle^{\frac12}}{\langle\xi_1
\rangle^{\frac{s}{2}}\langle\xi_3
\rangle^{\frac{s}{2}}} \, \, \, \mbox{ for all }\, \, \,
i=1,2,3.$$
Then by \eqref{se} and \eqref{mf2} we obtain
\begin{eqnarray*}
& &\int_*
\frac{\langle\xi_4\rangle^{s}\langle\xi_2\rangle^{1-s}}
{(\langle\xi_{12}\rangle
\langle\xi_{14}\rangle)^{1-b'}
\langle\xi_1\rangle^{s}\langle\xi_3\rangle^{s}
\prod_{j=1}^3\langle\tau_j+(-1)^j\xi_j^2
\rangle^{\frac12+}}
\prod_{j=1}^4F_j(\tau_j,\xi_j)\\
&\lesssim&\|\widetilde{{\cal F}}^{-1}(F_{4})\|_{L^{2}_{xt}}
\|\widetilde{{\cal F}}^{-1}\left(
\frac{\langle\xi\rangle^{\frac12}}
{\langle\tau+\xi^2
\rangle^{\frac12+}}F_{2}\right)\|_{L^{\infty}_{x}L^{2}_{t}}
\|\widetilde{{\cal F}}^{-1}\left(
\frac{\langle\xi\rangle^{-\frac{s}{2}}}{\langle\tau-\xi^2
\rangle^{\frac12+}}F_{3}\right)\|_{L^{4}_{x}L^{\infty}_{t}}\\
&\times&
\|\widetilde{{\cal F}}^{-1}\left(
\frac{\langle\xi\rangle^{-\frac{s}{2}}}{\langle\tau-\xi^2
\rangle^{\frac12+}}F_{1}\right)\|_{L^{4}_{x}L^{\infty}_{t}}
\lesssim
\prod_{j=1}^4\|F_j\|_{L^2}.
\end{eqnarray*}
\item
{\bf Case b):} In this case we borrow a power $\alpha=b'-\frac12+$
from the large denominator and we reduce our estimate to
\begin{eqnarray*}
\int_*
\frac{\langle\xi_4\rangle^{s}\langle\xi_2\rangle^{1-s}}
{\langle\xi_1\rangle^{s}\langle\xi_3\rangle^{s}
\prod_{j=1}^4\langle\tau_j+(-1)^j\xi_j^2
\rangle^{\frac12+}}
\prod_{j=1}^4F_j(\tau_j,\xi_j)\lesssim
\prod_{j=1}^4\|F_j\|_{L^2}.
\end{eqnarray*}
Again by symmetry we can assume that $|\xi_{1}|\geq |\xi_{3}|$.
We first observe that if the exponent of $\langle\xi_4\rangle$
were $\frac12$, then we could
simply use \eqref{se} for the function $F_{2}$ and \eqref{mf1}
for the function $F_{4}$ to obtain the estimate as we
did above. But in our case $s>\frac12$, so we have to do a bit more
work. We subdivide the analysis into subcases
\begin{itemize}
     \item
     {\bf Subcase 1):} $|\xi_{4}|\lesssim |\xi_{2}|$. In this case
we can write
$$\langle\xi_4\rangle^{s}\langle\xi_2\rangle^{1-s}\lesssim
\langle\xi_4\rangle^{\frac12}\langle\xi_2\rangle^{\frac12}$$
and we can indeed use \eqref{se} and \eqref{mf1}.
\item
{\bf Subcase 2):} $|\xi_{2}|\ll |\xi_{4}|$. Because we assumed that
$|\xi_{3}|\leq |\xi_{1}|$ and we are on the set
$\xi_{1}+\ldots+\xi_{4}=0$, it follows that
$|\xi_{4}|\lesssim|\xi_{1}| $. Then the estimate becomes
\begin{eqnarray*}
& &\int_*
\frac{\langle\xi_2\rangle^{1-s}}
{\langle\xi_3\rangle^{s}
\prod_{j=1}^4\langle\tau_j+(-1)^j\xi_j^2
\rangle^{\frac12+}}
\prod_{j=1}^4F_j(\tau_j,\xi_j)\\
&\lesssim&\|\widetilde{{\cal F}}^{-1}\left(
\frac{1}
{\langle\tau+\xi^2
\rangle^{\frac12+}}F_{4}\right)\|_{L^{4}_{xt}}
\|\widetilde{{\cal F}}^{-1}\left(
\frac{1}
{\langle\tau-\xi^2
\rangle^{\frac12+}}F_{1}\right)\|_{L^{4}_{xt}}
\|\widetilde{{\cal F}}^{-1}\left(
\frac{\langle\xi\rangle^{1-s}}{\langle\tau+\xi^2
\rangle^{\frac12+}}F_{2}\right)\|_{L^{\infty}_{x}L^{2}_{t}}\\
&\times&
\|\widetilde{{\cal F}}^{-1}\left(
\frac{\langle\xi\rangle^{-s}}{\langle\tau-\xi^2
\rangle^{\frac12+}}F_{3}\right)\|_{L^{2}_{x}L^{\infty}_{t}}
\lesssim
\prod_{j=1}^4\|F_j\|_{L^2},
\end{eqnarray*}
thanks   to \eqref{strichartz-p} for $p=2$, \eqref{se} and
\eqref{mf1}.
\end{itemize}
\end{itemize}
\end{proof}

\section{Proof of Proposition \ref{global}}\label{global-sec}
Based on Lemma \ref{dxe}, Theorem \ref{lwp} and  the arguments
presented
in \cite{ckstt:0}, Section 6 (see also the comments in
\cite{ckstt:0}, Section 7), the only result  that one needs to obtain is the
following
\begin{lemma}
     For any Schwartz function $w$, we have
\begin{equation}\label{inc-est}
\left|\int_T^{T+\delta} \Lambda_n(M_n; w(t))\ dt\right|
\lesssim \frac{1}{N^{2-}} \| Iw\|_{X^{1,3/4-}([T,T+\delta]\times \dbR)}^n
\end{equation}
for $n=6,8$, where $M_6$, $M_8$ are defined in Proposition
\ref{energy-increment}.
\end{lemma}
In \cite{ckstt:0} we were only able to obtain a decay of $N^{-1+}$,
which is why we could only prove global well-posedness
for $s>\frac23$.

The proof of this lemma is a corollary of the four lemmas that follow
in this section.

\begin{lemma}[$n=8$]\label{lem:M_8}
\begin{eqnarray*}
|M_8(\xi_1,\xi_2,\cdots,\xi_8)|\lesssim N_{soprano}m^{2}(N_{soprano}).
\end{eqnarray*}
\end{lemma}
This is a simple consequence of Lemma \ref{estm4}.
We now turn to the estimate of $\frac{d}{dt}E^{2}(Iw)$ involving
$\Lambda_{8}$.

\begin{lemma}
\begin{eqnarray*}
\left|\int_{T}^{T+\delta}\int\Lambda_8(M_8(\xi_1,\xi_2,\cdots,\xi_8))
\,dt\right|
\lesssim\frac{1}{N^{2-}}\|Iw\|_{1,\frac12+}^{8}.
\end{eqnarray*}
\end{lemma}
\begin{proof}
As in the proof of Lemma \ref{dxe}, also in this case
we first perform a Littlewood-Paley decomposition of the eight
factors $w$ so that the $\xi_{i}$ essentially are the constants
$N_{i}, \, i=1,\ldots,8$.
To recover the sum at the end we borrow a $N^{-\epsilon}_{soprano}$
from the large denominator $N_{soprano}$. Often this will not
be mentioned and it will only be recorded at the end by paying
a price equivalent to $N^{0+}$. Below we often use
the set of indices $R=\{soprano,alto,tenor\}$. Again we proceed by
analyzing  different cases:
\begin{itemize}
\item
{\bf Case a)} $N_{soprano}\sim N_{tenor}$.
By Lemma \ref{lem:M_8} and the fact that
$m(\xi)\langle\xi\rangle^{\frac12}$ is
increasing, we have
\begin{eqnarray*}
& &\left|\int_{T}^{T+\delta}\int\Lambda_8(M_8(\xi_1,\xi_2,\cdots,\xi_8))
\,dt\right|\\
&\lesssim&\sum_{R}\sum_{j_{1},\ldots,j_{5}\notin R}
\frac{N_{soprano}^{-2}}{m(N_{soprano})}\|D_{x}Iw_{soprano}\|_{L^6}
\|D_{x}Iw_{alto}\|_{L^6}
\|D_{x}Iw_{tenor}\|_{L^6}\\
& &\|D_{x}Iw_{j_{1}}\|_{L^6}
\Pi_{i=2,\ldots,5}\|D_{x}^{1/2-}Iw_{j_{i}}\|_{L^{\infty}}^2\lesssim
\frac{1}{N^{2-}}\|Iw\|_{1,\frac12+}^8.
\end{eqnarray*}
\item
{\bf Case b)} $N_{soprano}\gg N_{tenor}$.
By Lemma \ref{lem:M_8} and Lemma \ref{improved-strichartz}, and again the
monotonicity of $m(\xi)\langle\xi\rangle^{1/2}$, we have
\begin{eqnarray*}
\left|\int_{T}^{T+\delta}\int\Lambda_8(M_8(\xi_1,\xi_2,\cdots,\xi_8))
\,dt\right|&\lesssim&
N_{soprano}\|Iw_{soprano}w_{tenor}\|_{L^2}\|Iw_{alto}w_{baritone}
\|_{L^2}\\
&\times&\|w\|_{L^{\infty}}^4\lesssim\frac{1}{N^{2-}}
\|Iw\|_{1,\frac12+}^8.
\end{eqnarray*}
\end{itemize}
\end{proof}

\begin{lemma}[n=6]\label{lem:M_6}

     \noindent
\begin{itemize}
\item
If $N_{tenor}\gtrsim N$, we have
\begin{eqnarray}\label{eq:M_61}
|M_6(\xi_1,\xi_2,\cdots,\xi_6)|\lesssim
m(N_{soprano})^2N_{soprano}^2.
\end{eqnarray}
\item
If $N_{tenor}\ll N$, we have
\begin{eqnarray}\label{eq:M_62}
|M_6(\xi_1,\xi_2,\cdots,\xi_6)|\lesssim N_{soprano}N_{tenor}.
\end{eqnarray}
\end{itemize}
\end{lemma}
\begin{proof}
If $N_{soprano}\ll N$, $M_6$ vanishes.
Then we may assume $N_{soprano}\gtrsim N$.
Also in  the set $\xi_{1}+\ldots+\xi_{6}=0$ we have
$N_{alto}\sim N_{soprano}$.

The proof of (\ref{eq:M_61}) follows from \eqref{m4est}.
The proof of (\ref{eq:M_62}) is more delicate. By symmetry we
assume $soprano=1,~N_1\ge N_3\ge N_5,~N_2\ge N_4\ge
N_6$. Again we  analyze different cases.
\begin{itemize}
\item
{\bf Case a):} $alto=2$.
The MVT shows
\begin{eqnarray*}
|\sigma_6(\xi_1,\xi_2,\cdots,\xi_6)|\lesssim
m(N_1)^2N_1N_{12}+m(N_{tenor})^2N_{tenor}^2
\lesssim m(N_{soprano})^2N_{soprano}N_{tenor}.
\end{eqnarray*}
Next we estimate the second term in $M_6$
\begin{eqnarray*}
\sum(M_4(\xi_{abc},\xi_d,\xi_e,\xi_f)\xi_b+
M_4(\xi_a,\xi_{bcd},\xi_e,\xi_f)\xi_c+
M_4(\xi_a,\xi_b,\xi_{cde},\xi_f)\xi_d+
M_4(\xi_a,\xi_b,\xi_c,\xi_{def})\xi_e).
\end{eqnarray*}
Again by \eqref{m4est} one has that
\begin{equation}
|M_4(\xi_{abc},\xi_d,\xi_e,\xi_f)\xi_g|\lesssim
m(N_{soprano})^2N_{soprano}N_{tenor},
\label{good}\end{equation}
for every $a,\ldots,g\in \{1,\ldots,6\}$ and $g\ne soprano, alto$.
Thus we only have to consider the contributions
\begin{eqnarray*}
& &\left|\sum_{(a,e)\in \{3,5\}}\sum_{(d,f)\in \{4,6\}}
M_4(\xi_{a21},\xi_d,\xi_e,\xi_f)\xi_2+
M_4(\xi_{a},\xi_{21d},\xi_e,\xi_f)\xi_1\right|\\
&+&
\left|\sum_{(a,c)\in \{3,5\}}\sum_{(d,f)\in
\{4,6\}}M_4(\xi_{a},\xi_{12b},\xi_e,\xi_f)\xi_1+
M_4(\xi_a,\xi_{b},\xi_{12e},\xi_{f})\xi_2\right|\\
&+&
\left|\sum_{(a,c)\in \{3,5\}}\sum_{(d,f)\in
\{4,6\}}M_4(\xi_a,\xi_{b},\xi_{12c},\xi_f)\xi_2+
M_4(\xi_a,\xi_{b},\xi_{c},\xi_{12f})\xi_1\right|\\
&+&\left|\sum_{(a,e)\in \{3,5\}}\sum_{(d,f)\in \{4,6\}}
M_4(\xi_{a2c},\xi_d,\xi_{1},\xi_f)\xi_2+
M_4(\xi_{a},\xi_{2},\xi_c,\xi_{d1f})\xi_1\right|=\sum_{i=1}^{4}I_{i}.
\end{eqnarray*}
Observe first that  all the variables appearing  in the function
$M_{4}$ in $\sum_{i=1}^{3}I_{i}$ are strictly smaller that
$\frac{N}{2}$, hence
by \eqref{m41} it follows that
$$
\sum_{i=1}^{3}I_{i}\lesssim N_{soprano}N_{tenor}.
$$
To estimate $I_{4}$ we use \eqref{13m4} and the symmetry of
$M_{4}$. Then also in this case we obtain
$$I_{4}\lesssim N_{soprano}N_{tenor}.$$

\noindent
\item
{\bf Case b):} $alto=3$.
In this case we need some cancellation between the large terms coming
from $\sigma_{6}(\xi_{1},\ldots,\xi_{6})$ and the large terms
of the sum of the $M_{4}$. From \eqref{good} it is easy to see that
one needs to estimate only
\begin{eqnarray*}
  \widetilde{M_{6}}(\xi_{1},\ldots,\xi_{6})&=&-\frac{1}{6}(m_{1}^{2}
  \xi_{1}^{2}+m_{3}^{2}\xi_{3}^{2})\\
  &-&\frac{\xi_{1}}{36}\left(\sum_{(b,d,f)\in \{2,4,6\}}
  M_4(\xi_{a},\xi_{b1d},\xi_{3},\xi_f)+
  M_4(\xi_{a},\xi_b,\xi_{3},\xi_{d1f})\right)\\
  &-&
  \frac{\xi_{3}}{36}\left(\sum_{(b,d,f)\in \{2,4,6\}}
  M_4(\xi_{a},\xi_{b},\xi_{1},\xi_{d3f})+
  M_4(\xi_{a},\xi_{b3d},\xi_{1},\xi_{f})\right).
     \end{eqnarray*}
We now use \eqref{m4app} and the symmetries of $M_{4}$ to write
\begin{eqnarray*}
\widetilde{M_{6}}(\xi_{1},\ldots,\xi_{6})&=&-\frac{1}{6}(m_{1}^{2}
  \xi_{1}^{2}+m_{3}^{2}\xi_{3}^{2})\\
&-&\frac{\xi_{1}}{72}\left(\sum_{(b,d,f)\in \{2,4,6\}}
  \frac{m_{3}^{2}(\xi_{b1d}^{2}+\xi_{b1f}^{2})}{\xi_{3}}\right)
  +O(N_{soprano}N_{tenor})\\
  &-&\frac{\xi_{3}}{72}\left(\sum_{(b,d,f)\in \{2,4,6\}}
  \frac{m_{1}^{2}(\xi_{d3f}^{2}+\xi_{b3d}^{2})}{\xi_{1}}\right)
   +O(N_{soprano}N_{tenor})\\
&=&-\frac{1}{6}(m_{1}^{2}
  \xi_{1}^{2}+m_{3}^{2}\xi_{3}^{2})\\
&+&\frac{1}{72}\left(\sum_{(b,d,f)\in \{2,4,6\}}
  m_{3}^{2}(\xi_{b1d}^2+\xi_{b1f}^2)\right)
  +O(N_{soprano}N_{tenor})\\
  &+&\frac{1}{72}\left(\sum_{(b,d,f)\in \{2,4,6\}}
  m_{1}^{2}(\xi_{d3f}^{2}+\xi_{b3d}^{2})\right)
   +O(N_{soprano}N_{tenor})\\
  &=&-\frac{1}{72}m_{3}^{2}\sum_{(b,d,f)\in \{2,4,6\}}
  (\xi_{3}^{2}-\xi_{1bd}^{2})+(\xi_{3}^{2}-\xi_{1fb}^{2})\\
&-&\frac{1}{72}m_{1}^{2}\sum_{(b,d,f)\in \{2,4,6\}}
  (\xi_{1}^{2}-\xi_{3bf}^{2})+(\xi_{1}^{2}-\xi_{b3d}^{2})
  +O(N_{soprano}N_{tenor}),
     \end{eqnarray*}
and now it is clear that also in this case
$$|\widetilde{M_{6}}(\xi_{1},\ldots,\xi_{6})|\lesssim
N_{soprano}N_{tenor}.$$
\end{itemize}
\end{proof}
\begin{lemma}
\begin{eqnarray}\label{eq:A_6}
\left|\int_{T}^{T+\delta}\int\Lambda_6(M_6(\xi_1,\xi_2,\cdots,
\xi_6))\,dt\right|\lesssim\frac{1}{N^{2-}}
\|Iw\|_{1,\frac34-}^{6}.
\end{eqnarray}
\end{lemma}
\begin{proof}
Also in this case one uses a Littlewood-Paley decomposition to start.
We divide the proof into three different cases: {\bf Case a)} when
$N_{baritone}\gtrsim N$, {\bf Case b)}
when $N_{soprano}\ge N_{tenor}\gtrsim N\gg N_{baritone}$ and {\bf
Case c)}
when $N_{soprano}\sim N_{alto}\gtrsim N \gg N_{tenor}$. Below
we often use the two sets of indices
$S=\{soprano, alto,tenor, baritone\}$ and $R=\{soprano, alto,tenor\}$.
We also recall that thanks
to the fact that  $ m(\xi)|\xi|^{\frac12}$ is not decreasing,
\begin{equation}
     m(\xi)(1+|\xi|)\gtrsim \left\{ \begin{array}{l}
N, \, \, \, \mbox{ if } \, \, |\xi|>\frac{N}{2}\\
1, \, \, \, \mbox{ if } \, \, |\xi|\leq \frac{N}{2}.
\end{array}\right.
\label{m}\end{equation}

\begin{itemize}
\item
{\bf Case a):} $N_{baritone}\gtrsim N$.
By Lemma \ref{lem:M_6}, \eqref{m} and the Strichartz estimate
\eqref{strichartz-6}, we have
\begin{eqnarray*}
&&\left|\int_{T}^{T+\delta}\int\Lambda_6(M_6(\xi_1,\xi_2,\cdots,
\xi_6))\,dt\right|\lesssim\sum_{S}\sum_{j,k\notin S}
\frac{1}{m(N_{tenor})N_{tenor}^{1-}
m(N_{baritone})N_{baritone}^{1-}}\\
&\times&m(N_{soprano})N_{soprano}
\|w_{soprano}\|_{L^6}
m(N_{alto})N_{alto}\|w_{alto}\|_{L^6}\\
&\times&m(N_{tenor})N_{tenor}^{1-}\|w_{tenor}
\|_{L^6}m(N_{baritone})N_{baritone}^{1-}
\|w_{baritone}\|_{L^6}
\|Iw_j\|_{L^6}\|Iw_k\|_{L^6}\\
&\lesssim&\frac{1}{N^{2-}}
\sum_{N_{alto}\sim N_{soprano}}
\|Iw_{soprano}\|_{1,\frac12+}
\|Iw_{alto}\|_{1,\frac12+}\|Iw\|_{1,\frac12+}^{4}.
\end{eqnarray*}
and Cauchy-Schwarz with respect to $N_{alto}\sim N_{soprano}$ concludes 
the proof of this part.
\item
{\bf Case b):} $N_{soprano}\ge N_{tenor}\gtrsim N\gg N_{baritone}$.
This is the only part in which we need to use the space $X^{1,b}$
with $b\sim \frac34-$. By Lemma \ref{lem:M_6} and \eqref{m} we have
\begin{eqnarray*}
\left|\int_{T}^{T+\delta}\int\Lambda_6(M_6(\xi_1,\xi_2,\cdots,
\xi_6))\,dt\right|&\lesssim&\sum_{R}\sum_{j,h,k\in R^{c}}
\frac{1}{N_{tenor}m(N_{tenor})}
m(N_{soprano})N_{soprano}\|w_{soprano}w_{baritone}
\|_{L^2}\\
&\times&m(N_{alto})N_{alto}\|w_{alto}\|_{L^6}
m(N_{tenor})N_{tenor}\|w_{tenor}\|_{L^6}\\
&\times&\|D_{x}^{\frac12}Iw_j\|_{L^{12}}\|D_{x}^{\frac12}Iw_k\|_{L^{12}}
\|D_{x}^{\frac12}Iw_h\|_{L^{12}}.
\end{eqnarray*}
Using Lemma \ref{improved-strichartz}
and \eqref{m}, it is easy to
see that
$$m(N_{soprano})N_{soprano}\|w_{soprano}w_{baritone}
\|_{L^2}\lesssim N_{soprano}^{-1/2}\|Iw_{soprano}\|_{X^{1,\frac12+}}
\|Iw_{baritone}\|_{X^{1,\frac12+}}.$$
Also by the Sobolev inequalities and again \eqref{m}, 
$$\|D_{x}^{\frac12}Iw_j\|_{L^{12}}
\lesssim \|Iw_j\|_{X^{1,\frac12+}}, $$
and similarly for $h$ and $k$.
Collecting the above estimates one obtains
\begin{eqnarray*}
\left|\int_{T}^{T+\delta}\int\Lambda_6(M_6(\xi_1,\xi_2,\cdots,
\xi_6))\,dt\right|&\lesssim&
\frac{1}{N^{\frac32-}}\|Iw\|_{1,\frac12+}^{6}.
\end{eqnarray*}
Unfortunately the decay $N^{-\frac32+}$ is not enough for our
purposes. Because the local estimate allow us to handle
terms of  type $\|Iw\|_{1,\frac34-}$ (see Section 5), we take
advantage of the
extra denominators. To see this we use the identity
$$ \xi_{1}+\ldots+\xi_{4}=0 \, \, \Longrightarrow \, \,
\xi_{1}^{2}-\xi_{2}^{2}+\xi_{3}^{2}-\xi_{4}^{2}=2\xi_{12}\xi_{14},$$
proved in \cite{ckstt:0}. We consider only the case
  $N_{1}=N_{soprano}, N_{2}=N_{alto}$ and
$N_{3}=N_{tenor}$. Indeed if  $N_{5}=N_{tenor}$  the argument is
easier. Then in  the set $\xi_{1}+\ldots+\xi_{6}=0$
we  write
\begin{eqnarray*}
\sum_{i=1}^{6}(-1)^{i-1}\xi_{i}^{2}&=&
\xi_{1}^{2}-\xi_{2}^{2}+\xi_{3}^{2}-(\xi_{4}+\xi_{5}+\xi_{6})^{2}\\
& &+(\xi_{4}+\xi_{5}+\xi_{6})^{2}-\xi_{4}^{2}+\xi_{5}^{2}-\xi_{6}^{2}
\\
&=&2\xi_{12}\xi_{1456}+(\xi_{4}+\xi_{5}+\xi_{6})^{2}
-\xi_{4}^{2}+\xi_{5}^{2}-\xi_{6}^{2},
\end{eqnarray*}
which implies that
$$|\sum_{i=1}^{6}(-1)^{i-1}\xi_{i}^{2}|\gtrsim N^{2},$$
and for  $\lambda_{1}+\ldots+\lambda_{6}=0$
\begin{equation}
\label{crucial}    N^{2}\lesssim \max_{i=1,\ldots,6}
     |\lambda_{i}+(-1)^{i}\xi_{i}^{2}|.
\end{equation}
If the integral in time were performed on the whole real line
instead of $[T,T+\delta]$, then, after paying the price of the
extra factor
$\max_{i=1,\ldots,6}|\lambda_{i}+(-1)^{i}\xi_{i}^{2}|^{\frac14}$,
one would obtain
\begin{eqnarray*}
\left|\int_{T}^{T+\delta}\int\Lambda_6(M_6(\xi_1,\xi_2,\cdots,
\xi_6))\,dt\right|&\lesssim&
\frac{1}{N^{2-}}\|Iw\|_{1,\frac34-}^{6}.
\end{eqnarray*}
This argument has to be modified when the time integral is performed
on a finite interval $[T,T+\delta]$, due to the fact that
$\chi_{[T,T+\delta]}$, the characteristic function of the interval
$[T,T+\delta]$ is not smooth enough. A similar difficulty was
encountered also in \cite{ckstt:0}. We split
$$\chi_{[T,T+\delta]}(t)=a(t)+b(t),$$
where
$$\hat{a}(\tau)=\widehat{\chi_{[T,T+\delta]}}(\tau)
\eta(\tau/N^{2}),$$
and $\eta$ is supported on a small interval of 0 and equals
1 near 0, so $a$ is smoothing out $\chi_{[T,T+\delta]}$ at scale
$N^{-2}$. If one replaces $\chi_{[T,T+\delta]}(t)$ with
$a(t)$, then the argument above works because the Fourier transform
of $a(t)$ is supported on $|\tau|\ll N^{2}$ and one can still  obtain
the crucial inequality \eqref{crucial}. We now have to deal with
$b(t)$. It is easy to check that
$$\|b(t)\|_{L^{1}_{t}}\lesssim N^{-2}.$$
So we just have to show that
\begin{equation}
\label{supt}\sup_{t}|\Lambda_{6}(M_{6};w_{1}(t),\ldots,w_{6}(t))|
     \lesssim \prod_{j=1}^{6}\|Iw_j\|_{X^{1,\frac34-}}.
\end{equation}
We can crudely use Lemma \ref{lem:M_6} and obtain
\begin{eqnarray*}
     |\Lambda_{6}(M_{6};w_{1}(t),\ldots,w_{6}(t))|
     &\lesssim&
     m_{soprano}^{2}N_{soprano}^{2}
\|w_{soprano}\|_{L^{\infty}_{t}L^{2}_{x}}
\|w_{alto}\|_{L^{\infty}_{t}L^{2}_{x}}\\
&\times&\|w_{tenor}\|_{L^{\infty}_{t}L^{\infty}_{x}}
\|w_{baritone}\|_{L^{\infty}_{t}L^{\infty}_{x}}
     \prod_{j\notin S}\|Iw_j\|_{L^{\infty}_{t}L^{\infty}_{x}},
\end{eqnarray*}
which gives \eqref{supt} by the Sobolev embedding theorem.
\item
{\bf Case c):} $N_{soprano}\sim N_{alto}\gtrsim N \gg N_{tenor}$.
By Lemma \ref{lem:M_6}, Lemma \ref{improved-strichartz}, Sobolev
inequality and
\eqref{m}, we have
\begin{eqnarray*}
\left|\int\!\!\int\Lambda_6(M_6(\xi_1,\xi_2,\cdots,
\xi_6))\right|&\lesssim&\sum_{S}\sum_{j,h\notin S}
\frac{1}{m_{alto}^{2}N_{alto}}
N_{soprano}N_{tenor}\|Iw_{soprano}Iw_{tenor}
\|_{L^2}\\
&\times&N_{alto}\|Iw_{alto}Iw_{baritone}\|_{L^2}\|w_j\|_{L^{\infty}}
\|w_h\|_{L^{\infty}}
\lesssim \frac{1}{N^{2-}}\|Iw\|_{1,\frac12+}.
\end{eqnarray*}
\end{itemize}
This concludes the proof of the lemma.
\end{proof}

\end{document}